\documentclass[12pt,leqno]{article}
\usepackage{amsfonts}
\usepackage{amsmath, amsthm, amsfonts, amssymb, color}
\usepackage{mathrsfs}
\usepackage{color}
\theoremstyle{definition}

\newcommand{\scr}[1]{\mathscr #1}
\definecolor{wco}{rgb}{0.5,0.2,0.3}

\numberwithin{equation}{section} \theoremstyle{remark}

\newcommand{\ua}{\uparrow}

\title{{\bf Harnack Inequalities for   Stochastic Equations Driven by L\'evy Noise}\footnote{Supported in
 part by  Lab. Math. Com. Sys., NNSFC(11131003 and 11201073), SRFDP,
 the Fundamental Research Funds for the Central Universities, and the Program for New Century Excellent Talents in Universities
of Fujian (No.\ JA11051 and JA12053).} }
\author{
{\bf   Feng-Yu Wang$^{a), c)}$  and Jian Wang$^{b)}$}\\
\footnotesize{$^{a)}$School of Mathematical Sciences,
Beijing Normal
University, Beijing 100875, China}\\
\footnotesize{$^{b)}$School of Mathematics and Computer Science, Fujian Normal University, Fuzhou 350007, China}\\
 \footnotesize{$^{c)}$Department of Mathematics,
Swansea University, Singleton Park, SA2 8PP, United Kingdom}\\ \footnotesize{wangfy@bnu.edu.cn, F.-Y.Wang@swansea.ac.uk, jianwang@fjnu.edu.cn}}
\begin{document}
\allowdisplaybreaks
\def\R{\mathbb R}  \def\ff{\frac} \def\ss{\sqrt}
\def\N{\mathbb N} \def\kk{\kappa}
\def\dd{\delta} \def\DD{\Delta} \def\vv{\varepsilon} \def\rr{\rho}
\def\<{\langle} \def\>{\rangle} \def\GG{\Gamma} \def\gg{\gamma}
  \def\nn{\nabla} \def\pp{\partial} \def\E{\mathbb E}
\def\d{\text{\rm{d}}} \def\bb{\beta} \def\aa{\alpha} \def\D{\scr D}
  \def\si{\sigma} \def\ess{\text{\rm{ess}}}
\def\beg{\begin} \def\beq{\begin{equation}}  \def\F{\scr F}
\def\Ric{\text{\rm{Ric}}} \def\Hess{\text{\rm{Hess}}}
\def\e{\text{\rm{e}}} \def\ua{\underline a} \def\OO{\Omega}  \def\oo{\omega}
 \def\tt{\tilde} \def\Ric{\text{\rm{Ric}}}
\def\cut{\text{\rm{cut}}} \def\P{\mathbb P} \def\ifn{I_n(f^{\bigotimes n})}
\def\C{\scr C}      \def\aaa{\mathbf{r}}     \def\r{r}
\def\gap{\text{\rm{gap}}} \def\prr{\pi_{{\bf m},\varrho}}  \def\r{\mathbf r}
\def\Z{\mathbb Z} \def\vrr{\varrho} \def\ll{\lambda}
\def\L{\scr L}\def\Tt{\tt} \def\TT{\tt}\def\II{\mathbb I}
\def\i{{\rm in}}\def\Sect{{\rm Sect}}  \def\H{\mathbb H}
\def\M{\scr M}\def\Q{\mathbb Q} \def\texto{\text{o}} \def\LL{\Lambda}
\def\Rank{{\rm Rank}} \def\B{\scr B} \def\i{{\rm i}} \def\HR{\hat{\R}^d}
\def\to{\rightarrow}\def\l{\ell}\def\1{\mathbf 1}

\maketitle

\begin{abstract}By using coupling argument and regularization approximations of the
underlying subordinator, dimension-free Harnack inequalities are
established for a class of stochastic equations driven by a L\'evy
noise containing a subordinate Brownian motion. The Harnack
inequalities    are
 new even for  linear equations driven by L\'evy noise, and the gradient estimate implied by
our log-Harnack inequality considerably generalizes some recent
results on gradient estimates and coupling properties derived for
L\'evy processes or linear equations driven by L\'evy noise. The
main results are also extended to semi-linear stochastic equations
in Hilbert spaces.
\end{abstract} \noindent
 \vskip 2cm

\section{Introduction}

Due to their broad range of applications in heat kernel estimates, functional inequalities, transportation-cost inequalities and properties of invariant measures, the dimension-free Harnack inequality with powers introduced in
\cite{W97} and the  log-Harnack inequality introduced in
\cite{RW} have been intensively investigated for
stochastic (partial) differential equations driven by Gaussian
noises, see   \cite{wangbook13,W11} and references therein.  However, due to technical
difficulty on construction of couplings for jump processes, the
study for stochastic equations driven by purely jump L\'evy noise is very limited. The only known results on this type of
Harnack inequalities are presented in \cite{W11,WJ} for linear
stochastic differential equations (i.e. O-U processes) driven by
purely jump L\'evy processes, where \cite{W11} uses   coupling
through the Mecke formula and \cite{WJ} adopts known heat kernel
bounds of the $\aa$-stable processes. Recently, using regularization
approximations of the time-change, Zhang established in \cite{Zhang}
the Bismut formula for stochastic differential equations with
Lipschitz continuous drifts driven by the $\aa$-stable process. In
this paper, we will make use of Zhang's argument together with a
coupling method to derive Harnack inequalities for   stochastic
equations driven by L\'evy noise, which provide explicit heat kernel estimates (see Remark 2.2 below).

Let $W:=(W_t)_{t\ge 0}, S:=(S(t))_{t\ge 0}$ and $V:=(V_t)_{t\ge 0}$ be independent stochastic processes, where $W$ is the Brownian motion on $\R^d$ with $W_0=0$; $V$ is a locally bounded
measurable process on $\R^d$ with $V_0=0$; and $S$ is the subordinator induced by a Bernstein function $B$, i.e., $S$ is a  one-dimensional non-negative increasing L\'evy process with $S(0)=0$ and
$\E \e^{-r S(t)}= \e^{-t B(r)},\, t,r\ge 0.$ Then $(W_{S(t)})_{t\ge 0}$ is a L\'evy process with symbol $\psi:=B(|\cdot|^2).$

We consider the following stochastic equation on $\R^d$:
\beq\label{1.1}  X_t= X_0+ \int_0^t b_s(X_s)\d s+\int_0^t \si_s \d
W_{S(s)} +V_t,\ \ t\ge 0,\end{equation} where $\si:
[0,\infty)\to\R^d\otimes\R^d$ is measurable and locally bounded, and
$b: [0,\infty)\times \R^d\to \R^d$ is measurable, locally bounded
and continuous in the second variable.

We shall need the following conditions on
$\si$ and $b$: \beg{enumerate}
\item[{\bf (H1)}] $\si_t^{-1}$ exists and is locally bounded, i.e.
there exists an increasing function $\ll$ on $[0,\infty)$ such that
$\|\si_t^{-1}\|\le \ll_t$ for $t\ge 0.$
\item[{\bf (H2)}] There exists a locally bounded measurable function $K$ on $[0,\infty)$ such that
\beq\label{H*} \<b_t(x)-b_t(y), x-y\>\le K_t |x-y|^2,\ \ x,y\in\R^d, t\ge 0.\end{equation} \end{enumerate}

It is easy to see that {\bf (H2)} implies the existence, uniqueness and  non-explosion of the solution, see e.g. the proof of  \cite[Theorem 177]{ST}. Now, for any $x\in \R^d$, let $(X_t(x))_{t\ge 0}$ be the unique
solution to (\ref{1.1}) for $X_0=x.$ We aim to establish Harnack
inequalities for the associated Markov operator $P_t$ on
$\B_b(\R^d)$:
$$P_t f(x):=\E f(X_t(x)),\ \ t\ge 0, x\in \R^d, f\in \B_b(\R^d).$$

Comparing with the O-U type equations studied in \cite{W11, WJ}, our equation is with a more general time-dependent drift. Moreover, to compare the L\'evy term in (\ref{1.1}) with those in \cite{W11,WJ} under a lower bound condition of the L\'evy measure, we may replace $W_{S(t)}$ in (\ref{1.1}) by a L\'evy process $L_t$ with  L\'evy measure  $\nu(\d x)\ge c |x|^{-d} B(|x|^{-2})\d x$ for some constant $c>0.$ In fact,  in this case we may split $L_t$ into two independent L\'evy parts, where one of them has L\'evy measure $c |x|^{-d} B(|x|^{-2})\d x$ and is thus a subordinate Brownian motion (cf.\ \cite{W,BRW}), and the integral of $\si$ w.r.t. the other  can be combined with the term $V_t$.

In Section 2 we state our main results, which are then proved in Section 3 by using regularization approximations of $S(t)$ and the coupling by change of measure. Finally,  the main results are extended in Section 4 to semilinear SPDEs by using finite-dimensional approximations.

\section{Main Results}

\beg{thm}\label{T1.1} Assume {\bf (H1)} and {\bf (H2)}, and let
$K(t)= \int_0^s K_u\d u$ for $t\ge 0.$  \beg{enumerate}\item[$(1)$]  For any $T>0$ and strictly positive $f\in \B_b(\R^d),$
\beg{equation*} \beg{split}P_T\log f(y)&\le \log P_T f(x) +\ff{|x-y|^2}2\inf_{t\in (0,T]} \E \bigg\{\ff{\ll_t^2}{\int_0^t\e^{-2K(s)}\d S(s)}\bigg\}, \ \ x,y\in \R^d.\end{split}\end{equation*}
\item[$(2)$] For any $T>0$ and $f\in \B_b(\R^d)$,
$$|\nn P_T f|(x)^2 \le \{P_Tf^2(x)-(P_Tf(x))^2\} \inf_{t\in [0,T]} \E \bigg\{\ff{\ll_t^2}{\int_0^t\e^{-2K(s)}\d S(s)}\bigg\},\ \ x\in \R^d,$$ where $|\nn P_T f|(x)$ is the local Lipschitz constant of $P_T f$ at
point $x$, i.e.
$$|\nn P_T f|(x)=\limsup_{y\to x} \frac{|P_Tf(y)-P_Tf(x)|}{|y-x|}.$$
\item[$(3)$]  For any  $T>0, p>1$ and positive $f\in \B_b(\R^d),$
\beg{equation*}\beg{split} (P_T f(y))^p&\le (P_T f^p(x))\bigg(\E\inf_{t\in (0,T]} \exp\bigg[\ff{p\ll_t^2|x-y|^2}{2(p-1)^2 \int_0^t\e^{-2K(s)}\d S(s)}\bigg]\bigg)^{p-1},\ \ x,y\in\R^d.\end{split}\end{equation*}
\end{enumerate}\end{thm}

\paragraph{Remark 2.1.} (1) When $S(t)=t$ and $V_t=0$, the equation (\ref{1.1}) reduces to the SDE driven by Brownian motion. In this case the assertions in Theorem \ref{T1.1} coincide with the corresponding ones derived in the diffusion setting, e.g. when $\si_t=\sqrt{2}I$ and $b_t=b$, assertions (1), (2) and (3) reduce respectively to (1.3), (iii) and (1.2) in \cite{W10} for $M=\R^d$. These inequalities are sharp as they are equivalent to the underlying curvature condition, see \cite[Theorem 1.1]{W10}.

(2) For general subordinator $S$, Theorem \ref{T1.1} improves
\cite[Theorem 1.1]{RW03} and \cite[Theorems 1.2 and 1.3]{ORW} for
generalized Mehler semigroups.   Theorem \ref{T1.1} (1) and (3) are new even for linear stochastic equations driven by L\'evy
processes, for which the Harnack inequality has been investigated in
\cite{W11} by using the Mecke formula. Since in \cite{W11} the
density of the L\'evy measure was used, so that  the derived inequalities
can not be extended to infinite-dimensions as we did in Section 4
for our present Harnack inequalities.

(3) The gradient
inequality in Theorem \ref{T1.1}(2) generalizes the main results in
\cite{W,BRW,SSW} for L\'evy processes or linear equations driven by
L\'evy noise. When $K\le 0$ and $\ll$ is bounded,  it follows from
Theorem \ref{T1.1}(2) that $$|P_T f(x)-P_Tf(y)|\le \|\ll\|_\infty
\|f\|_\infty |x-y| \ss{\E\ff 1 {S(T)}},\ \ f\in\B_b(\R^d), f\ge 0.$$
This implies the coupling property provided $\E \ff 1 {S(T)}\to 0$
as $T\to\infty.$ Thus, the main results in \cite{W11',BRW, SW, SW2}
on the coupling property of L\'evy processes or linear equations
driven by L\'evy noise are generalized, see also \cite{WJ12} for the
recent study of the coupling property of L\'{e}vy processes with
drift. If furthermore $K\le -\theta$ for some constant $\theta>0$,
we have $|X_t(x)-X_t(y)|\le\e^{-\theta t}|x-y|,$ which together with
Theorem \ref{T1.1}(2), implies the exponential convergence of $\nn
P_T$: there exists a constant $C>0$ such that (see the proof of
\cite[Theorem 1.1]{W} for details)
$$|\nn P_T f|(x)^2\le C \e^{-2\theta T} \{P_Tf^2(x)- (P_Tf(x))^2\},\ \
x\in\R^d, T\ge 1, f\in \B_b(\R^d).$$

To illustrate Theorem \ref{T1.1}, we consider the equation driven by stable like processes. In the following result $W_{S(t)}$ is the $\aa$-stable process when $B(r)=r^{\aa/2},\aa\in (0,2).$

\beg{cor}\label{C1.2} Assume that the assumptions of Theorem $\ref{T1.1}$ hold and assume that there exist some constants $\theta\in (0,1)$ and $c,r_0>0 $ such that $B(r)\ge c r^\theta$ holds for any $r\ge r_0$. Then there exists a constant $C>0$ such that
\beg{enumerate}\item[$(1)$]  For any $T>0$ and strictly positive $f\in \B_b(\R^d),$
$$P_T\log f(y)\le \log P_T f(x) +\ff{C|x-y|^2} {(T\land 1)^{\ff 1 \theta}},\ \ x,y\in \R^d;$$
\item[$(2)$] For any $T>0$ and $f\in\B_b(\R^d)$,
$$|\nn P_Tf|(x)^2\le\{P_Tf^2(x)-(P_Tf(x))^2\} \ff{C}{(T\land 1)^{\ff 1 \theta}}, \ \ x\in \R^d;$$
\item[$(3)$]  When $\theta\in (\ff 1 2,1)$,  any  $T>0, p>1, x,y\in\R^d$ and positive $f\in \B_b(\R^d),$
$$ (P_T f(y))^p
\le (P_T f^p(x))\exp\Bigg[\ff{Cp|x-y|^2}{(p-1)(T\land 1)^{\ff 1 \theta}}+\ff{C(p|x-y|^2)^{\ff{\theta}{2\theta-1}}}{[(p-1)(T\land 1)]^{\ff 1 {2\theta-1}}}\Bigg]. $$
\end{enumerate}
\end{cor}

\paragraph{Remark 2.2} Among some other applications of Harnack inequalities summarized in \cite[\S 1.4]{wangbook13} (see also \cite{W10, WY}), to save space we only mention here heat kernel estimates implied by Corollary \ref{C1.2}. Let $p_T(x,y)$ be the density of $P_T$ with respect to the Lebesgue measure (the existence is well known as the equation is non-degenerate). By \cite[Proposition 3.1(4)]{WY} for $\mu$ replacing by the Lebesgue measure (note that the proof works also for $\si$-finite quasi-invariant measures), Corollary \ref{C1.2} (1) and (3) imply the entropy inequality
\beq\label{HE1} \int_{\R^d} p_T(x,z)\log\ff{p_T(x,z)}{p_T(y,z)}\d z\le \ff{C|x-y|^2}{(T\land 1)^{\ff 1 \theta}},\ \ x,y\in \R^d, T>0,\end{equation} and when $\theta\in (\ff 1 2,1)$, for any $p>1$ and $T>0$,
\beq\label{HE2} \int_{\R^d} p_T(x,z)\Big(\ff{p_T(x,z)}{p_T(y,z)}\Big)^{\ff 1 p}\d z\le \exp\bigg[\frac{C|x-y|^2}{(p-1)(T\land 1)^{\ff 1 \theta}}+ \ff{Cp^{\ff{1-\theta}{2\theta-1}}|x-y|^{\ff{2\theta}{2\theta-1}}}{[(p-1)(T\land 1)]^{\ff 1 {2\theta-1}}}\bigg],\ \ x,y\in\R^d.\end{equation}
Moreover, it is obvious that Corollary \ref{C1.2}(2) implies
\beq\label{HE3} \sup_{x\in\R^d} \int_{\R^d} |\nn\log p_T(\cdot,y)(x)|^2p_T(x,y)\d y\le \ff {C}{(T\land 1)^{\ff 1 \theta}},\ \ T>0.\end{equation}

\section{Proofs of Theorem \ref{T1.1} and Corollary \ref{C1.2}}

We first explain the main idea of the proof. As in \cite{Zhang} we consider the following regularization of $S$:
$$S_\vv(t):=\ff 1 \vv \int_t^{t+\vv} S(s)\d s +\vv t,\ \ t\ge 0, \vv>0.$$ Then $S_\vv$ is strictly increasing, absolutely continuous and $S_\vv\downarrow S$ as $\vv\downarrow 0.$ For each $\vv>0$, we consider the approximation equation
\beq\label{2.0} X_t^\vv =X_0+\int_0^t b_s(X_s^\vv)\d s +\int_0^t\si_s\d W_{S_\vv(s)} +V_t,\ \ t\ge 0.\end{equation} Since $S_\vv$ is absolutely continuous, this equation is indeed driven by the Brownian motion so that we are able to establish the Harnack inequalities for the associated operator $P_t^\vv$ by using coupling. Finally, by proving $P_t^\vv\to P_t$ as $\vv\to 0$, we derive the corresponding Harnack inequalities for $P_t$.

\subsection{The case of absolutely continuous time-change}
Let $\ell$ be an absolutely continuous and strictly increasing function on $[0,\infty)$ with $\ell(0)=0,$ and let $v: [0,\infty)\to \R^d$ be measurable and  locally bounded  with $v_0=0$.
We consider the following equation
\beq\label{2.1} X_t^{\ell,v}   =X_0+\int_0^t b_s(X_s^{\ell,v})\d s +\int_0^t\si_s\d W_{\ell(s)} +v_t,\ \ t\ge 0.\end{equation}
Under our general assumptions, this equation has a unique solution. Let
$$P_t^{\ell,v} f(x)= \E f(X_t^{\ell,v}(x)),\ \ x\in\R^d, t\ge 0, f\in\B_b(\R^d),$$ where $X_t^{\ell,v}(x)$ is the solution to (\ref{2.1}) for $X_0=x.$

Now, for fixed  $T>0$ and $x,y\in\R^d$, we intend to construct a coupling to derive   the Harnack inequalities of $P_T^{\ell,v}$. To this end,
 let $(Y_t)_{t\ge 0}$ solve the equation
\beq\label{2.2} Y_t= y+ \int_0^t b_s(Y_s)\d s +\int_0^t \si_s \d W_{\ell(s)} +v_t +\int_0^t  \big\{\xi_s1_{[0,\tau)}(s)\big\}\cdot \ff{X^{\ell,v}_s-Y_s}{|X^{\ell,v}_s-Y_s|} \d\ell(s),\end{equation} where
$$\xi_t:=\ff{|x-y|\exp[-\int_0^t K_s\d s]}{\int_0^T \exp[-2\int_0^tK_s \d s]\d\ell(t)},\ \ t\ge 0$$ and
$$\tau:=\inf\{t\ge 0: X^{\ell,v}_t=Y_t\}.$$ To construct a solution to (\ref{2.2}),  we consider the equation
\beq\label{2.2'} \tt Y_t= y+ \int_0^t b_s(\tt Y_s)\d s +\int_0^t \si_s \d W_{\ell(s)} +v_t +\int_0^t  \big\{\xi_s1_{\{X^{\ell,v}_s\ne\tt Y_s\}}\big\}\cdot \ff{X^{\ell,v}_s-\tt Y_s}{|X^{\ell,v}_s-\tt Y_s|} \d\ell(s).\end{equation} Since $(z,z')\mapsto \ff{z-z'}{|z-z'|}$ is locally Lipschitz continuous on the domain $\{(z,z')\in \R^d\times\R^d: z\ne z'\}$, the joint equation of (\ref{2.1}) and (\ref{2.2'}) has a unique solution up to the coupling time
$$ \tt\tau:=\inf\{t\ge 0: X^{\ell,v}_t=\tt Y_t\}.$$ Let
$$Y_t=  \tt Y_t1_{[0, \tt\tau)}(t)+
X^{\ell,v}_t1_{[\tt\tau,\infty)}(t). $$ Then $(Y_t)_{t\ge 0}$ is a solution to (\ref{2.2}) with $\tau=\tt\tau$.

\beg{lem} \label{L2.1}For the above constructed coupling $(X^{\ell,v}_t,Y_t)_{t\ge 0}$, there holds $\tau\le T$, i.e. $X^{\ell,v}_T=Y_T.$\end{lem}

\beg{proof} By (\ref{2.1}) and (\ref{2.2}) we have
$$\d(X^{\ell,v}_t-Y_t)= (b_t(X^{\ell,v}_t)-b_t(Y_t))\d t -\xi_t\cdot \ff{X^{\ell,v}_t-Y_t}{|X^{\ell,v}_t-Y_t|} \d\ell(t),\ \ t<\tau.$$ Then {\bf (H2)} and the absolutely continuity of $\ell$ yield
$$\d |X^{\ell,v}_t-Y_t|^2 \le 2 K_t |X^{\ell,v}_t-Y_t|^2\d t -2 \xi_t |X^{\ell,v}_t-Y_t|\d\ell(t),\ \  t<\tau.$$ Note that for two continuous semimartingales $M_t$ and $\tt M_t$, the inequality $\d M_t\le \d\tt M_t$ means that they have the same martingale part and $\tt M_t-M_t$ is an increasing process.  Thus,
$$\d\big\{|X^{\ell,v}_t-Y_t|\e^{-\int_0^t K_s\d s}\big\}\le -\xi_t \e^{-\int_0^t K_s\d s}\d\ell(t),\ t<\tau.$$ Therefore, if $\tau>T$ then
$$0<|X^{\ell,v}_T-Y_T|\e^{-\int_0^T K_s\d s}\le |x-y|-\int_0^T \xi_t \e^{-\int_0^t K_s\d s}\d t =0,$$ which is a contradiction. \end{proof}

To derive the Harnack inequality, we define
$$\tt W_t:= W_t+ \int_0^{t\land \ell(\tau)} \ff{\xi_{\ell^{-1}(s)}}{|X^{\ell,v}_{\ell^{-1}(s)}-Y_{\ell^{-1}(s)}|}   \si^{-1}_{\ell^{-1}(s)} (X^{\ell,v}_{\ell^{-1}(s)} -Y_{\ell^{-1}(s)}) \d s,\ \ t\ge 0.$$
By the Girsanov theorem, $(\tt W_t)_{t\ge 0}$ is the $d$-dimensional Brownian motion under the probability $\d\Q:=R\d\P$, where
\beg{equation*}\beg{split} &R:=\exp\bigg[-\int_0^{\ell(\tau)}\<\eta_t,\d W_t\>-\ff 1 2 \int_0^{\ell(\tau)}|\eta_t|^2\d t\bigg],\\
& \eta_t:= \ff{\xi_{\ell^{-1}(t)}}{|X^{\ell,v}_{\ell^{-1}(t)}-Y_{\ell^{-1}(t)}|}   \si^{-1}_{\ell^{-1}(t)} (X^{\ell,v}_{\ell^{-1}(t)} -Y_{\ell^{-1}(t)}),\ \ t\in [0,\ell(\tau)).\end{split}\end{equation*}
Reformulating (\ref{2.2}) by
$$Y_t=y+\int_0^tb_s(Y_s)\d s + \int_0^t\si_s\d\tt W_{\ell(s)} +v_t,\ \ t\ge 0,$$  we conclude from the definition of $P_t^{\ell,v}$ and Lemma \ref{L2.1} that
\beq\label{2.3} P_T^{\ell,v} f(y)=\E_\Q f(Y_T) = \E[Rf(Y_T)]= \E[R f(X^{\ell,v}_T)].\end{equation}
It is now more or less standard that  this formula implies the following result.

\beg{prp}\label{P2.2} For any strictly positive $f\in \B_b(\R^d)$,
\beq\label{LHL} P_T^{\ell,v} \log f(y)\le \log P_T^{\ell,v}f(x) +\inf_{t\in (0,T]}\ff{\ll_t^2 |x-y|^2}{2\int_0^t \e^{-2K(s)}\d \ell(s)},\ \ x,y\in\R^d,\end{equation} and for any $p>1$,
\beq\label{Hl} (P_T^{\ell,v} f)^p(y)\le (P_T^{\ell,v}f^p(x)) \inf_{t\in (0,T]}\exp\bigg[\ff{p\ll_t^2|x-y|^2}{2(p-1) \int_0^t \e^{-2K(s)}\d \ell(s)}\bigg],\ \ x,y\in\R^d.\end{equation} \end{prp}

\beg{proof} (1) By (\ref{2.3}) and the Young inequality that for any probability measure $\nu$ on $\R^d$, if  $g_1, g_2\ge 0$ with $\nu(g_1) = 1$, then  $$\nu(g_1g_2) \le  \log \nu(e^{g_2})+\nu(g_1 \log g_1),$$ we obtain
$$P_T^{\ell,v}\log f(y)=\E[R\log f(X^{\ell,v}_T)]\le \log P_T^{\ell,v} f(x) +\E[R\log R].$$ By the definitions of $R,\eta_t,\xi_t$ and noting that $\tau\le T$, we have
\beg{equation*}\beg{split} \E[R\log R] &= \E_\Q\log R=\E_\Q\bigg\{ -\int_0^{\ell(\tau)} \<\eta_t,\d \tt W_t\>+\ff 1 2 \int_0^{\ell(\tau)} |\eta_t|^2\d t\bigg\}\\
&= \ff 1 2 \E\int_0^{\ell(\tau)} |\eta_t|^2\d t \le \ff{\ll_T^2} 2 \int_0^{\ell(T)} \xi_{\ell^{-1}(t)}^2\d t
 =\ff{\ll_T^2} 2 \int_0^T \xi_t^2\d\ell(t)\\
 & =\ff{\ll_T^2|x-y|^2}{2\int_0^T \exp[-2\int_0^tK_s\d s]\d\ell(t)}.\end{split}\end{equation*}
This implies that
$$P_T^{\ell,v} \log f(y)\le \log P_T^{\ell,v}f(x) + \ff{\ll_T^2 |x-y|^2}{2\int_0^T\e^{-2K(s)}\d \ell(s)},\ \ x,y\in\R^d.$$
Now, for any $t\in (0,T]$, let $P^{\ell,v}_{t,T} f(x)= \E (f(X_{t,T}^{\ell,v}(x))$, where $(X_{t,T}^{\ell,v}(x))_{t\le T}$ solves the equation
$$X_{t,T}^{\ell,v}(x)= x +\int_t^Tb_s(X_{t,s}^{\ell,v}(x))\d s +\int_t^T \si_s\d W_{\ell(s)} + v_T-v_t,\ \ t\le T.$$  By the Markov property we have $P^{\ell,v}_T= P^{\ell,v}_tP^{\ell,v}_{t,T}$. So, applying the above inequality to $t$ and $P^{\ell,v}_{t,T}f$ in place of $T$ and $f$ respectively, and noting that by the Jensen inequality
$$ P_t^{\ell,v} \log P^{\ell,v}_{t,T} f \ge P_t^{\ell,v} P_{t,T}^{\ell,v} \log f=P_T^{\ell,v}\log f,$$ we obtain
 (\ref{LHL}).

(2) By (\ref{2.3}) and the H\"older inequality, we obtain
$$(P_T^{\ell,v}f(y))^p= (\E[Rf(X^{\ell,v}_T)])^p\le (P_T^{\ell,v}f^p(x)) (\E R^{\ff p{p-1}})^{p-1}.$$ Since
\beg{equation*}\beg{split} \E R^{\ff p{p-1}}&= \E\bigg(\exp\bigg[-\ff p{p-1}\int_0^{\ell(\tau)} \<\eta_t,\d W_t\> -\ff p{2(p-1)}\int_0^{\ell(\tau)} |\eta_t|^2\d t\bigg]\bigg)\\
&= \E\bigg(\exp\bigg[-\ff p{p-1}\int_0^{\ell(\tau)} \<\eta_t,\d W_t\> -\ff {p^2}{2(p-1)^2}\int_0^{\ell(\tau)} |\eta_t|^2\d t\bigg]\\
&\qquad\,\,\,\times\exp\bigg[ \ff p{2(p-1)^2}\int_0^{\ell(\tau)} |\eta_t|^2\d t \bigg]\bigg)\\
&\le \exp \bigg[\ff{p\ll_T^2}{2(p-1)^2} \int_0^{\ell(T)} \xi_{\ell^{-1}(t)}^2 \d t\bigg]= \exp\bigg[\ff{p\ll_T^2 |x-y|^2}{2(p-1)^2 \int_0^T\e^{-2K(s)}\d\ell(s)}\bigg],\end{split} \end{equation*} we obtain
$$(P_T^{\ell,v} f)^p(y)\le (P_T^{\ell,v}f^p(x)) \exp\bigg[\ff{p\ll_T^2 |x-y|^2}{2(p-1) \int_0^T\e^{-2K(s)}\d \ell(s)}\bigg],\ \ x,y\in \R^d.$$ From this we obtain (\ref{Hl}) similarly as in the first part of the proof.
\end{proof}

\subsection{Proofs of Theorem \ref{T1.1} and Corollary \ref{C1.2}} To prove Theorem \ref{T1.1} using Proposition \ref{P2.2}, we need the following lemma to ensure that $X_t^{(n)}\to X_t$ as $n\to \infty$, where $X^{(n)}$ solves (\ref{2.0}) for $\vv= \ff 1 n.$

\beg{lem}\label{L2.3} Assume that
\beg{enumerate} \item[$(i)$] $\si$ is piecewise  constant: there exists a sequence $\{t_n\}_{n\ge 0}$ with $t_0=0$ and $t_n\uparrow\infty$  such that $\si_\cdot= \sum_{i=1}^\infty  1_{[t_{i-1},t_i)} \si_{t_{i-1}};$
\item[$(ii)$] $b$ is globally Lipschitzian: for any $T>0$ there exists a constant $C>0$ such that
$$|b_t(x)-b_t(y)|\le C|x-y|,\ \ t\in [0,T], x,y\in\R^d.$$\end{enumerate} Then $\lim_{n\to \infty}X_t^{(n)}=X_t$ holds for all $t>0.$\end{lem}

\beg{proof} Let $T>0$ be fixed. By $(i)$ and $(ii)$, for any $t\in [0,T]$, we have
\beq\label{GG}\beg{split}  |X_t^{(n)}-X_t|\le& C\int_0^t |X_s^{(n)}-X_s|\d s \\
&+2\sup_{t\in [0,T]} \|\si_t\| \Big(\sum_{  t_i<t}  |W_{S_{1/n}(t_i)}-W_{S(t_i)}| +|W_{S_{1/n}(t)}-W_{S(t)}|\Big).\end{split}\end{equation}
Moreover, it is easy to see that $(ii)$ and the local boundedness of $b,\si, V$ imply
$ \sup_{n\ge 1} |X_t^{(n)}|<\infty$ and thus,
$\phi_t:=\limsup_{n\to\infty} |X_t^{(n)}-X_t|<\infty$ for any $t\in [0,T].$ Combining these with (\ref{GG}) and using the Fatou lemma and the fact that $S_{1/n}\downarrow S$ as $n\uparrow\infty$, we arrive at
$\phi_t\le C\int_0^t \phi_s\d s$ for any $t\in [0,T].$ Therefore, $\phi_t=0$ for all $t\in [0,T]$ and the proof is thus finished. \end{proof}

\beg{proof}[Proof of Theorem \ref{T1.1}] According to
\cite[Proposition 2.3]{ATW12}, (2) follows from (1). So, we only
prove (1) and (3). To apply Lemma \ref{L2.3}, we shall make
approximations of $b$ and $\si$. Since $C_b(\R^d)$ is dense in
$L^1(P_T(x,\cdot)+P_T(y,\cdot))$, where $P_T(z,\d z')$ is the
transition probability for $P_T$, we only consider strictly positive
$f\in C_b(R^d)$.

(a) We first assume that   $(i)$ and $(ii)$ in Lemma \ref{L2.3} hold. By applying Proposition \ref{P2.2} to $P_T^{S_{1/n}, V}$ and noting that Lemma \ref{L2.3} implies
$$P_Tf =\lim_{n\to\infty} \E P_T^{S_{1/n},V}f,\ \ f\in C_b(\R^d),$$ we obtain (\ref{LHL}) and (\ref{Hl}) for $(S,V)$ in place of $(\ell,v).$
Then the log-Harnack inequality follows by taking expectations to (\ref{LHL}), and the Harnack inequality with power follows by taking expectations to (\ref{Hl}) and using the H\"older inequality:
\beg{equation*}\beg{split} P_T f(y) &= \E P_T^{S,V}f(y)\le \E \bigg\{(P_T^{S,V}f^p(x))^{\ff 1 p} \inf_{t\in [0,T]} \exp\Big[\ff{\ll_t^2 |x-y|^2}{2(p-1)\int_0^t \e^{-2K(s)}\d S(s)}\Big]\bigg\}\\
&\le \big(\E P_T^{S,V} f^p(x)\big)^{\ff 1 p} \bigg(\E \inf_{t\in [0,T]} \exp\Big[\ff{p\ll_t^2|x-y|^2}{2(p-1)^2\int_0^t \e^{-2K(s)}\d S(s)}\Big]\bigg)^{\ff{p-1}p}\\
&= \big(  P_T f^p(x)\big)^{\ff 1 p} \bigg(\E \inf_{t\in [0,T]} \exp\Big[\ff{p\ll_t^2|x-y|^2}{2(p-1)^2\int_0^t \e^{-2K(s)}\d S(s)}\Big]\bigg)^{\ff{p-1}p}.\end{split}\end{equation*}

(b) Assume that $(ii)$ in Lemma \ref{L2.3} holds. Since for a fixed
sample of $S$, the class of piecewise constant functions on
$[0,\infty)$ is dense in $L^2_{loc}(\d S)$, we may find out a
sequence of $\R^d\otimes\R^d$-valued functions $\{\si^{(n)}\}_{n\ge
1}$ satisfying $(i)$ in Lemma \ref{L2.3} such that $\si^{(n)}\to
\si$ in $L^2([0,T];\d S)$ and $\|(\si_t^{(n)})^{-1}\|\le \ll_T$ for
$t\in [0,T].$ Let $\tt X_t^{(n)}$ solve (\ref{1.1}) for $\si^{(n)}$
in place of $\si$, and let $\tt P_t^{(n)}$ be the associated Markov
operator. According to (a), the assertions in Theorem \ref{T1.1} hold
for $\tt P_T^{(n)}$ in place of $P_T$. By $(ii)$ we have
\beq\label{G*0} |X_t-\tt X_t^{(n)}|\le C\int_0^t|X_s-\tt
X_s^{(n)}|\d s + \bigg|\int_0^t (\si_s^{(n)}-\si_s)\d
W_{S(s)}\bigg|,\  \ t\in [0,T].\end{equation} Since $\si^{(n)}\to
\si$ in $L^2([0,T];\d S),$ we have (see e.g. \cite[Theorem 88(v) on page
53]{ST})
$$\lim_{n\to\infty}  \E^S \bigg|\int_0^t (\si_s^{(n)}-\si_s)\d W_{S(s)}\bigg|^2=\lim_{n\to\infty} \int_0^t \|\si_s^{(n)}-\si_s\|_{HS}^2\d S(s) =0,$$ where $\E^S$ is the conditional expectation given $S$. Then, as in the proof of Lemma \ref{L2.3}, by letting $n\to\infty$ in (\ref{G*0}) we obtain
$\lim_{n\to\infty} \E^S |\tt X_T^{(n)}-X_T|=0,$ so that
$$P_T f= \E P_T^{S,V} f = \lim_{n\to\infty} \E\big(\E^{S,V} f(\tt X_T^{(n)})\big) =\lim_{n\to\infty}\tt P_T^{(n)} f.$$ Therefore, Theorem \ref{T1.1} also holds for $P_T$.

(c) In general, let $\tt b_t(x)= b_t(x)-K_tx.$ Then (\ref{H*}) is equivalent to the dissipative property of $\tt b$:
$$\<\tt b_t(x)-\tt b_t(y), x-y\>\le 0,\ \ t\ge 0, x,y\in\R^d.$$ Let $(\tt b^{(n)})_{n\ge 1}$ be the Yoshida approximation of $\tt b$, i.e.
$$ \tt b_t^{(n)}(x)= n\Big\{\Big(I-\ff 1 n \tt b_t\Big)^{-1}(x)-x\Big\},\ \ t\ge 0, x\in \R^d.$$ Then (see e.g.\ \cite[Section 2]{DRW09}), $\tt b^{(n)}$ is dissipative and globally Lipschitzian in the sense of $(ii)$ of Lemma \ref{L2.3}, $|\tt b^{(n)}|\le |\tt b|$, and  $\lim_{n\to\infty} \tt b^{(n)}=\tt b.$ Let $b^{(n)}_t(x)=\tt b^{(n)}_t(x) + K_t x.$ Then, $b_t^{(n)}$ satisfies $(ii)$, and (\ref{H*}) holds for $b_t^{(n)}$ in place of $b_t$.

Now, let $\bar X_t^{(n)}$ solve (\ref{1.1}) for $b^{(n)}$ in place of $b$. Then, according to (b), the associated Markov operator $\bar P_t^{(n)}$ satisfies the claimed inequalities in Theorem \ref{T1.1}.  So, if
\beq\label{**W} \lim_{n\to\infty} \bar P_T^{(n)} f= P_Tf,\ \ f\in C_b(\R^d),\end{equation} then  we complete the proof by applying Theorem \ref{T1.1} to $\bar P_T^{(n)}$ and letting $n\to\infty$. The proof of (\ref{**W}) is straightforward by the constructions of $\tt b$ and $\tt b^{(n)}$, from which we have
\beq\label{*D}\beg{split} \d |X_t-\bar X_t^{(n)}|^2&= 2 K_t |X_t-\bar X_t^{(n)}|^2\d t + 2 \<\tt b_t(X_t)-\tt b_t^{(n)}(\bar X_t^{(n)}),X_t-\bar X_t^{(n)}\>\d t\\
& \le  2 K_t |X_t-\bar X_t^{(n)}|^2\d t + 2 \<\tt b_t(X_t)-\tt b_t^{(n)}(X_t),X_t-\bar X_t^{(n)}\>\d t\\
&\le (2 K_t+1) |X_t-\bar X_t^{(n)}|^2\d t + |\tt b_t(X_t)-\tt b_t^{(n)}(X_t)|^2\d t.\end{split}\end{equation} Let $\tau_m=\inf\{t\ge 0: |X_t|\ge m\}$ for $m\ge 1$. We obtain  from (\ref{*D}) that
 $$|X_{T\land\tau_m}-\bar X^{(n)}_{T\land\tau_m}|^2\le \e^{2\int_0^T(K_t+1)\d t} \int_0^{T\land\tau_m}|\tt b_t(X_t)-\tt b_t^{(n)}(X_t)|^2\d t.$$
Since $\{|\tt b_t(X_t)-\tt b^{(n)}_t(X_t)|: t\le T\land\tau_m, n\ge 1\}$ is bounded and $\lim_{n\to\infty} \tt b_t^{(n)}= \tt b_t$, this implies  $\lim_{n\to\infty} |X_{T\land\tau_m}-\bar X_{T\land\tau_m}^{(n)}|^2=0$  for all $m\ge 1$. Combining this with  $\tau_m\uparrow\infty$ as $m\uparrow\infty$, we conclude that $\lim_{n\to\infty}\bar X_T^{(n)}=X_T$ a.s. and thus prove (\ref{**W}).  \end{proof}

\beg{proof}[Proof of Corollary \ref{C1.2}]   By Theorem \ref{T1.1}, it suffices to prove for $T\in (0,1].$   There exists a constant $c_1\ge 1$ such that  for any $k\ge 1$,
\beg{equation}\label{W0*}\beg{split} \E \ff 1 {S(t)^k}&= \ff 1{\GG(k)} \int_0^\infty r^{k-1} \e^{-tB(r)}\d r \le \ff{\exp[cr_0^\theta t]}{\GG(k)} \int_0^\infty r^{k-1}\exp[{-cr^\theta t}]\d r\\
&= \ff{\exp[cr_0^\theta t]}{\GG(k)c^{\ff 1 \theta}}\int_0^\infty
r^{k-1} \e^{-t r^\theta}\d r \le  c_1\E\ff 1 {\tt S (t)^k},\ \ t\in (0,1],\end{split}\end{equation} where $ \tt
S$ is the subordinator associated to the Bernstein function
$r\mapsto r^\theta.$ Therefore, (1) and (2) follow  from Theorem
\ref{T1.1}(1)-(2) and \cite[(2.2)]{GRW} by noting that
\beq\label{*00} \ff 1{\int_0^t \e^{-2K(s)}\d S(s)}\le \ff{c_2}{S(t)},\ \ t\in (0,1]\end{equation} holds for some constant $c_2>0.$  To prove (3), we make use the third display from below  in the proof of \cite[Theorem 1.1]{GRW} for
$\kk=1$, i.e.
$$\E \e^{\ll/\tt S(t)} \le 1+\bigg(\exp\Big[\ff{c_3\ll^{\ff{\theta}{2\theta-1}}}{t^{\ff 1 {2\theta-1}}}\Big]-1\bigg)^{\ff{2\theta-1}\theta}\le
\exp\Big[ \ff{c_4\ll}{t^{\ff 1 \theta}}+\ff{c_4\ll^{\ff\theta{2\theta-1}}}{t^{\ff 1{2\theta-1}}}\Big],\ \ \ll,t\ge 0$$ for some constants $c_3,c_4>0.$ This along with (\ref{W0*}) yields that
\beg{equation*}\beg{split} \E \e^{\ll/ S(t)} &\le 1+ c_1\E\big(\e^{\ll/\tt S(t)}-1\big)\le \E\e^{c_1\ll/\tt S(t)}\\
&\le  \exp\Big[\ff{c_5\ll}{t^{\ff 1 \theta}}+\ff{c_5\ll^{\ff\theta{2\theta-1}}}{t^{\ff 1{2\theta-1}}}\Big]\end{split}\end{equation*} for some $c_5>0.$ Combining this with (\ref{*00}) we prove (3) from Theorem \ref{T1.1}(3).   \end{proof}

\section{Extension to semi-linear SPDEs}

Let $(H,\<\cdot,\cdot\>,|\cdot|)$ be a separable Hilbert space,
$V:=(V_t)_{t\ge 0}$ be a locally bounded measurable stochastic
process on $H$, $W=(W_t)_{t\ge 0}$ be a cylindrical Brownian motion
on $H$, and $S=(S_t)_{t\ge 0}$ be a one-dimensional non-negative
increasing L\'evy process associated to a Bernstein function $B$ as
introduced in Section 1. Recall that $W$ can be formally formulated
as \beq\label{QP} W_t= \sum_{i=1}^\infty B_t^i e_i,\end{equation}
where $\{B^i\}_{i\ge 1}$ is a family of independent one-dimensional
Brownian motions, and $\{e_i\}_{i\ge 1}$ is an orthonormal basis of
$H$. Thus, for any orthonormal family $\{e_i\}_{i=1}^n$, the process
$(\<W,e_1\>,\cdots, \<W,e_n\>)$ is a Brownian motion on $\R^n$. As
in the finite-dimensional case, we assume that $W, S$ and $V$ are
independent. Let $\scr L(H)$ be the set of all bounded linear
operators on $H$.

Consider the following stochastic equation on $H$:
\beq\label{3.1} X_t= \e^{A t} X_0+ \int_0^t \e^{A(t-s)}F_s(X_s)\d s +\int_0^t \e^{A(t-s)}\si_s \d W_{S(s)} +V_t,\ \ t\ge 0,\end{equation}
where $\si: [0,\infty)\to \scr L(H)$ is measurable and locally bounded, and $A$ and $F$ satisfy
\beg{enumerate} \item[{\bf (A1)}] $(A,\D(A))$ is a negatively definite self-adjoint operator on $H$ such that
\beq\label{*3} \int_0^t \|\e^{sA}\|_{HS}^2\d s <\infty,\ \  t>0;\end{equation}
\item[{\bf (A2)}] $F: [0,\infty)\times H\to H$ is measurable, bounded on bounded sets and satisfies
\beq\label{GL2} |F_s(x)-F_s(y)|\le K_s |x-y|,\ \ x,y\in H, s\ge 0\end{equation} holds for some locally bounded measurable function $K$ on $[0,\infty).$
  \end{enumerate}

Note that if $V_t=\int_0^t \e^{(t-s)A} \tt\si_s\d L_s$ holds for some $\tt \si: [0,\infty)\to \L(H)$ and some noise $L$, (\ref{3.1}) is   known as the definition of mild solutions to the stochastic differential equation
$$\d X_t= AX_t\d t + F_t(X_t)\d t + \si_t\d W_{S(t)} +\tt\si_t\d L_t.$$

We first confirm the existence and uniqueness of the solution
\eqref{3.1}. By {\bf (A1)}, the operator $-A$ has discrete spectrum
with eigenvalues $0\le\rho_1\le\rho_2\le\cdots\le \rho_n\uparrow
\infty$. From now on, we let $\{e_i\}_{i\ge 1}$ be the corresponding
eigenbasis, i.e. an orthonormal basis of $H$ such that $Ae_i=-\rho_i
e_i, i\ge 1.$

\beg{prp}\label{PP} Assume {\bf (A1)}-{\bf (A2)}.
Then\beg{enumerate}\item[$(1)$] Let $B^k:=\<W,e_k\>$ as in
(\ref{QP}). Then
$$Y_t:= \int_0^t \e^{(t-s)A}\si_s \d W_{S(s)} = \sum_{k,j=1}^\infty \bigg(\int_0^t\<\e^{(t-s)A}\si_se_k, e_j\>\d B^k_{S(s)}\bigg) e_j,\ \ t\ge 0,$$
gives rise to a stochastically continuous process on $H$ such that
$\E^S \int_0^T |Y_t|^2\d t<\infty$ for $T>0$, where $\E^S$ is the
conditional expectation given $S$. In particular, $Y$ has a
measurable modification.
\item[$(2)$] Fix a measurable modification of $Y:= (Y_{t})_{t\ge 0}$, and denote it again by $$Y_t=\int_0^t \e^{(t-s)A}\si_s \d W_{S(s)}, \ \ t\ge 0.$$  For any $X_0\in H$,
the equation $(\ref{3.1})$ has a unique solution.  \end{enumerate} \end{prp}

\beg{proof}   (1) We first prove  \beq\label{GO} U_t:=
\sum_{k,j=1}^\infty \int_0^t \<\e^{A(t-s)}\si_s e_k, e_j\>^2\d S(s)
<\infty \ \text{a.s.},\ \ t\ge 0,\end{equation} and $\int_0^T U_t\d
t<\infty$ for any $T>0.$ Note that for each $t\ge 0$, $U_t<\infty$
a.s.  implies that $Y_t$ is a well defined $H$-valued random
variable with $\E^S |Y_t|^2=U_t<\infty$ (see e.g. \cite[Theorem 88(v) on page
53]{ST}), and  $\int_0^T U_t\d
t<\infty$ implies $\E^S \int_0^T |Y_t|^2\d t<\infty.$

It is easy to see that
\beq\label{D0} \beg{split} U_t &= \int_0^t \|\e^{A(t-s)}\si_s\|_{HS}^2\d S(s)\\
&\le \Big(\sup_{s\in [0,t]}\|\si_s\|^2\Big) \int_0^t
\|\e^{A(t-s)}\|_{HS}^2\d S(s)=: C_t \gg_t,\ \ t\ge
0.\end{split}\end{equation} It then follows from (\ref{*3}) that
\beg{equation*}\beg{split}   \int_0^T  U_t\d t &\le C_T \int_0^T\gg_t\d t =C_T\int_0^T\d S(s)\int_s^T \|\e^{(t-s)A}\|_{HS}^2\d t  \\
&\le C_T  S(T)\int_0^T\|\e^{tA}\|_{HS}^2\d t
<\infty.\end{split}\end{equation*} So, $Y_\cdot\in L^2([0,T]\to H;\d
t)$ and $  U_t<\infty$ a.s.  for a.e.-$t\ge 0$. It remains to show
that $U_t<\infty$ a.s. for all $t\ge 0$, so that $Y_t$ is an
$H$-valued random variable for each $t\ge 0$. For any $t>0$, there
exists $t'\in (0,t)$ such that $\gg_{t'},\gg_{t-t'}<\infty$ a.s.
Since $S=  S(\cdot+t')-S(t')$ in law,  $\int_{t'}^{t}
\|\e^{A(t-s)}\|_{HS}^2\d S(s)=\gg_{t-t'}=\int_0^{t-t'}
\|\e^{A(t-t'-s)}\|_{HS}^2\d S(s)$ in law as well. Thus,
\beg{equation*}\beg{split} \gg_t  &= \int_0^{t'} \|\e^{A(t-s)}\|_{HS}^2 \d S(s) +\int_{t'}^{t} \|\e^{A(t-s)}\|_{HS}^2 \d S(s)\\
& = \gg_{t'} + \int_{t'}^{t} \|\e^{A(t-s)}\|_{HS}^2\d S(s)<\infty,\
\text{a.s.}\end{split}\end{equation*} Therefore, by (\ref{D0}) we
have $ U_t<\infty$ a.s. for all $t\ge 0.$

To prove the stochastic continuity of $Y$, we note that for any
$t\ge 0$ and $h>0$, we have
$$|Y_{t+h}-Y_t|\le |\e^{hA}Y_t-Y_t| +|I(h)|,$$ and
$$I(h):=\int_t^{t+h} \e^{(t+h-s)A} \si_s \d W_{S(s)}.$$
Note that $C_t= \sup_{s\in [0,t]}\|\si_s\|^2<\infty.$ We have, for
$h\in (0,1),$
$$ \E^S |I(h)|^2 \le C_{t+1} \int_{t-h}^t \|\e^{(t-s)A}\|_{HS}^2 \d
S(s+h)$$ which  in  law equals to  $C_{t+1} \int_{t-h}^t
\|\e^{(t-s)A}\|_{HS}^2 \d S(s).$ Since $U_t<\infty$ a.s. and
$S(t)=S(t-)$ a.s. for fixed $t$, we conclude that $I(h)\to 0$ in
probability as $h\to 0.$ Therefore, for any $\vv>0$,
$$\limsup_{h\downarrow 0} \P(|Y_{t+h}-Y_t|\ge \vv) \le
\limsup_{h\downarrow 0}\Big\{ \P\Big(|\e^{hA}Y_t- Y_t|\ge \ff
\vv 2\Big)+ \P\Big(|I(h)|\ge\ff \vv 2\Big)\Big\}=0.$$ Similarly, we
can prove $\lim_{s \uparrow t}\P(|Y_t-Y_s|\ge\vv)=0$ for any
$t,\vv>0.$ Due to the stochastic continuity, the process $Y$ has a
measurable modification (see \cite[Theorem 3]{Am}).

(2) Once a measurable modification of $Y$ is fixed, as explained in
Section 1, we let $\tt X_t= X_t-Y_t-V_t$ and reformulate (\ref{3.1})
as
 $$\tt X_t= \e^{At}X_0+\int_0^t \e^{A(t-s)}F_s(\tt X_s +Y_s+V_s)\d s,\ \ t\ge 0,$$ which has a unique solution due to {\bf (A2)}. \end{proof}

 We note that in the proof of Proposition \ref{PP}, for different measurable modifications of $Y$, the corresponding solutions derived for the equation (\ref{3.1}) are equivalent,
  i.e. they are modifications each other as well. When $V=0$ and $\si_s$ is independent of $s$ with $\si e_i=\bb_ie_i$ holding for some sequence $\{\bb_i\}\subset \R$,
   solutions to (\ref{3.1}) have been investigated in \cite{PZ}.

 By Proposition \ref{PP}, we define
$$P_t f(x)= \E f(X_t(x)),\ \ t\ge 0, x\in H, f\in \B_b(H),$$ where $X(x)$ is the solution to (\ref{3.1}) for $X_0=x.$
We shall make use of finite-dimensional approximations to derive the Harnack inequalities from Theorem \ref{T1.1}.
 Note that $P_t$ is independent of modifications of $Y$, and  is thus unique due to Proposition \ref{PP}.

For $n\ge 1$, let $H_n=\text{span}\{e_1,\cdots, e_n\}$, and let $\pi_n$ be the orthogonal projection from $H$ onto $H_n$. Let
$$A^{(n)}= \pi_n A, \ F^{(n)}= \pi_n F,\ \si^{(n)}= \pi_n\si,\ W^{(n)}= \pi_n W,\ V^{(n)}= \pi_n V.$$ For any $n\ge 1,$ consider the following equation on $H_n$:
\beq\label{4.2} X_t^{(n)} = \pi_n X_0 +\int_0^t \big\{A^{(n)}
X_s^{(n)}+ F_s^{(n)}(X_s^{(n)})\big\}\d s +\int_0^t \si_s^{(n)}\d
W_{S(s)}^{(n)} +V_t^{(n)},\ \ t\ge 0.\end{equation} Let $P_t^{(n)}$
be the associated Markov operator. It is easy to see that assertions
in Theorem \ref{T1.1} hold  for $P_T^{(n)}$. Letting $n\to\infty$,
we conclude that assertions in Theorem \ref{T1.1}, and hence in
Corollary \ref{C1.2}, hold for the present $P_T$.

\beg{thm} Assume {\bf (A1)} and {\bf (A2)}. If $\|(\si_t^{(n)})^{-1}\|\le \ll_t$ for some increasing function $\ll$ on $[0,\infty)$ and large $n$, then assertions in Theorem \ref{T1.1}  and   Corollary \ref{C1.2}  hold for $H$ in place of $\R^d$.
\end{thm}

\beg{proof} As explained above we only consider positive  $f\in
C_b(H)$. In this case, by the assertions for $P_T^{(n)}$ and the
dominated convergence theorem, it suffices to prove that $
\lim_{n\to\infty} X_T^{(n)} = X_T$  in  law. This follows from
$$ \lim_{n\to\infty}\E^{S,V} |X_t^{(n)}-X_t|=0,\ \
t>0,$$ which can be easily verified as in the proof of \cite[Theorem 2.1]{WZ}.
\end{proof}

Similarly to the finite-dimensional situation, if $P_T$ has a quasi-invariant measure $\mu$ then according to \cite[Proposition 3.1]{WY}, the assertions in Corollary \ref{C1.2} imply that $P_T$ has a heat kernel $p_T(x,y)$ with respect to $\mu$ and estimates (\ref{HE1}), \eqref{HE2} and \eqref{HE3} hold for $\mu$ and $H$ replacing the Lebesgue measure and $\R^d$ respectively.

\paragraph{Acknowledgement.} The authors would like to thank the associate editor,  the two referees, Dr.\ Jianhai Bao and Dr.\ Shaoqing Zhang for useful comments and corrections.

\end{document}